\documentclass[reqno,12pt]{amsart}

\usepackage{epsf}
\usepackage{graphics}
\usepackage{amssymb}
\usepackage{amsmath}

\date{}

\theoremstyle{plain}
\newtheorem{theorem}{Theorem}

\theoremstyle{definition}
\newtheorem*{definition}{Definition}

\theoremstyle{remark}

\def\R{{\mathbb R}}

\title{Legendrian framings for two-bridge links} 

\author{Sebastian Baader \and Masaharu Ishikawa}

\begin{document}

\begin{abstract} We define the Thurston-Bennequin polytope of a two-component link as the convex hull of all pairs of integers that arise as framings of a Legendrian representative. The main result of this paper is a description of the Thurston-Bennequin polytope for two-bridge links. As an application, we construct non-quasipositive surfaces in $\R^3$ all whose sub-annuli are quasipositive.
\end{abstract}

\maketitle

\section{Introduction}

A Legendrian knot in the standard contact space $\R^3$ defines an integer, called Thurston-Bennequin number, via its framing. The maximal Thurston-Bennequin number of a knot is obtained by maximizing this integer over all Legendrian representatives of that knot. Precise values of the maximal Thurston-Bennequin number are known for some special classes of knots, such as positive knots or two-bridge knots (\cite{Ta}, \cite{Ng}). Given a link $L \subset \R^3$ with two components, we may ask which pairs of integers can be realised as framings of a Legendrian representative. The convex hull of these pairs of integers defines a polytope in $\R^2$, which we call the Thurston-Bennequin polytope $\Delta(L)$ of the link $L$. As an example, the Thurston-Bennequin polytope of the negative torus link of type $T(2,-4)$ is shown in Figure~1 (this is a special case of Theorem~1 below).
\begin{figure}[ht]
\scalebox{1.0}{\raisebox{-0pt}{$\vcenter{\hbox{\epsffile{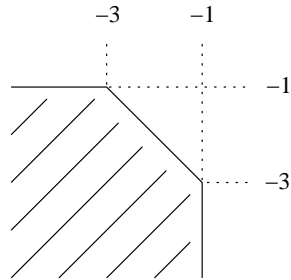}}}$}} 
\caption{$\Delta(T(2,-4))$}
\end{figure}
In this note we will determine the Thurston-Bennequin polytope for oriented two-bridge links. Our description makes use of the maximal Thurston-Bennequin number for oriented links, denoted by $\overline{tb}$, and the linking number $lk$. However, the polytope itself does not depend on orientation.

\begin{theorem} Let $L$ be a two-bridge link with two components $L_1$, $L_2$. The Thurston-Bennequin polytope of $L$ is
$$\Delta(L)=\\
\{(x_1,x_2) \in \R^2|x_1 \leq -1, \; x_2 \leq -1, \; x_1+x_2 \leq \overline{tb}(L)-2lk(L)\}.$$
Moreover, every integer point of $\Delta(L)$ is realised by a Legendrian framing.
\end{theorem}

Theorem~1 has an interesting application concerning quasipositive surfaces in $\R^3$. Quasipositive surfaces can be defined as the minimal family of embedded surfaces in $\R^3$ containing the positive Hopf band (i.e. the fibre of the positive Hopf link), which is closed under the plumbing operation and under taking incompressible sub-surfaces (\cite{Ru1}, \cite{Ru2}). Here a sub-surface is called incompressible if its embedding into the ambient surface induces an injective map on the level of fundamental groups. By definition, every incompressible sub-surface of a quasipositive surface is quasipositive. In particular, this is true for all incompressible sub-annuli. It is tempting to conjecture the converse: a surface embedded in $\R^3$ is quasipositive if all its incompressible sub-annuli are. In the case of pretzel surfaces, this is true \cite{Ru3}. In general, it is false.

\begin{theorem} There exist non-quasipositive connected surfaces in $\R^3$ all whose incompressible sub-annuli are quasipositive.
\end{theorem}

We will construct a surface with the above properties in the third section. The next section is devoted to the proof of Theorem~1, which is based on Ng's work \cite{Ng}.

\section{Construction of Legendrian representatives}

The standard contact structure on $\R^3$ is induced by the 1-form $dz-y \, dx$. As mentioned above, a Legendrian knot has a natural framing induced by the contact planes. The Thurston-Bennequin number is defined as the linking number of the link induced by that framing. In order to define the Thurston-Bennequin number for oriented links, it is convenient to introduce the front projection of links, i.e. the orthogonal projection to the $y$-$z$-plane. Front projections of Legendrian links contain cusps, as shown in Figure~2. The Thurston-Bennequin number $tb(\mathcal{L})$ of an oriented Legendrian link $\mathcal{L} \subset \R^3$ is defined as 
\begin{equation}
tb(\mathcal{L})=w(D)-\frac{1}{2} cu(D),
\label{tb_def}
\end{equation}
where $w(D)$ and $cu(D)$ denote the writhe (i.e. the algebraic crossing number) and the number of cusps of the front diagram $D$ of $\mathcal{L}$, respectively. For example, the Thurston-Bennequin number of the Legendrian link of Figure~2 is $-3$ or $-7$, depending on the choice of orientation. Let $\mathcal{L}=\mathcal{L}_1 \cup \mathcal{L}_2$ be a Legendrian representative of a link $L$ with two components $L_1$, $L_2$. The following equality is an immediate consequence of the above definition:
$$tb(\mathcal{L})=tb(\mathcal{L}_1)+tb(\mathcal{L}_2)+2lk(L).$$
Therefore, if a pair of integers $(x_1, x_2)$ is induced by a Legendrian framing of $L$, then the following three inequalities hold:
$$x_1 \leq \overline{tb}(L_1), \; x_2 \leq \overline{tb}(L_2), \; x_1+x_2 \leq \overline{tb}(L)-2lk(L).$$

\begin{proof}[Proof of Theorem~1]
Let $L$ be a two-bridge link with components $L_1$ and $L_2$. Since these components are unknotted, we have
$$\overline{tb}(L_1)=\overline{tb}(L_2)=-1.$$
Theorem~1 is settled once we have shown that every integer point on the diagonal line segment with endpoints $(-1,\overline{tb}(L)-2lk(L)+1)$, $(\overline{tb}(L)-2lk(L)+1,-1)$ is realised by a Legendrian framing. Indeed, if a pair of integers $(x_1, x_2)$ is realised, then so are the pairs $(x_1-1, x_2)$ and $(x_1, x_2-1)$, by simple stabilisations. 

Two-bridge links admit a description by finite continued fractions. In~\cite{Ng}, Ng proves that every two-bridge link has a diagram corresponding to an alternating continued fraction
$$[a_1, -a_2,\ldots,(-1)^{n-1}a_n]=a_1+\frac{1}{-a_2+\ldots+\frac{1}{(-1)^{n-1}a_n}}$$
with $a_i \geq 2$, for all $i \in \{1,\ldots,n\}$.
These diagrams in turn define a canonical front projection $F(a_1, \ldots, a_n)$, called Legendrian rational front, as illustrated in Figure~2.
\begin{figure}[ht]
\scalebox{1.2}{\raisebox{-0pt}{$\vcenter{\hbox{\epsffile{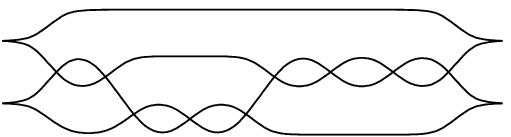}}}$}} 
\caption{$F(2,3,4)$}
\end{figure}
Let $D=F(a_1, \ldots, a_n)$ be a rational front projection for a Legendrian representative $\mathcal{L}$ of the link $L$. According to Ng, $D$ realises the maximal Thurston-Bennequin number of $L$. We observe that the upper component of $D$, say $\mathcal{L}_1$, has Thurston-Bennequin number $-1$, so the corresponding point of the Thurston-Bennequin polytope has coordinates $(-1, \overline{tb}(L)-2lk(L)+1)$. We will gradually move from this point along the diagonal line towards the point $(\overline{tb}(L)-2lk(L)+1,-1)$, via a suitable family of Legendrian front diagrams.

First we will arrange all the numbers $a_i$ corresponding to self-crossings of the lower component $\mathcal{L}_2$ to be even. For this purpose we use the fact that the value of a continued fraction is preserved if we replace a segment of the form $a_i, -a_{i+1}$ by the segment $a_i-1,1,a_{i+1}-1$ and change the signs of all subsequent coefficients, provided $a_i, a_{i+1} \geq 2$. This replacement can be done independently at all self-crossings (note that if $a_i$ corresponds to self-crossings of $\mathcal{L}_2$, then $a_{i+1}$ does not). However, this produces a non-alternating continued fraction: the crossings corresponding to coefficients $1$ carry the wrong sign. We may still associate a front diagram to this fraction by introducing a so-called zig-zag to the component $\mathcal{L}_2$ at all single crossings with the wrong sign. Figure~3 shows the two possible ways of introducing a zig-zag at a crossing.

\begin{figure}[ht]
\scalebox{1.0}{\raisebox{-0pt}{$\vcenter{\hbox{\epsffile{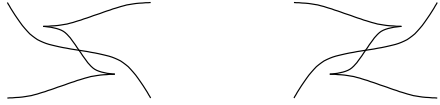}}}$}} 
\caption{two kinds of zig-zag}
\end{figure}

In this manner we obtain a front diagram $D'$ which still maximises the Thurston-Bennequin number of $L$. In fact all three numbers $tb(\mathcal{L})$, $tb(\mathcal{L}_1)$, $tb(\mathcal{L}_2)$ are the same for both diagrams $D$, $D'$. This follows from~(\ref{tb_def}) and the fact that all self-crossings of $\mathcal{L}_2$ are negative (the increase in the writhe from $D$ to $D'$ equals the number of cancelled self-crossing of $\mathcal{L}_2$ and is therefore neutralised by the additional pairs of cusps). 

Next we wish to modify $D'$ so that the number of self-crossings of $\mathcal{L}_2$ between the lower two strands is at least as high as the number of self-crossings of $\mathcal{L}_2$ between the middle two strands. If this is not yet the case, we may simply turn the diagram $D'$ upside down and drag the lower strand to the top, see Figure~4. Again, this operation does not affect the numbers $tb(\mathcal{L})$, $tb(\mathcal{L}_1)$, $tb(\mathcal{L}_2)$.

\begin{figure}[ht]
\scalebox{1.0}{\raisebox{-0pt}{$\vcenter{\hbox{\epsffile{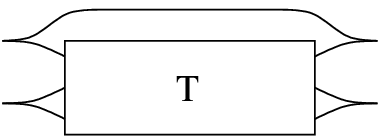}}}$}} \quad $\longrightarrow$ \quad
\scalebox{1.0}{\raisebox{-0pt}{$\vcenter{\hbox{\epsffile{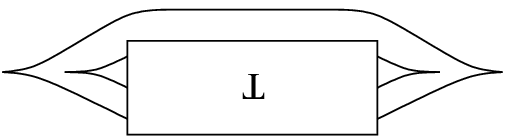}}}$}}
\caption{}
\end{figure}

At last, we will describe two kinds of operations on the diagram $D'$ that do change the numbers $tb(\mathcal{L}_1)$ and $tb(\mathcal{L}_2)$.
\begin{enumerate}
\item Replace a zig-zag on the component $\mathcal{L}_2$ by a zig-zag on the component $\mathcal{L}_1$, (compare Figure~3). This operation changes the Thurston-Bennequin numbers of $\mathcal{L}_1$ and $\mathcal{L}_2$ by $-1$ and $+1$, respectively.\\
\item Replace two consecutive self-crossings of $\mathcal{L}_2$ at the lower level by two consecutive self-crossings of $\mathcal{L}_1$ at the upper level, as shown in Figure~5. This operation corresponds to a double flype, but not to a Legendrian isotopy of front diagrams. In fact the Thurston-Bennequin numbers of $\mathcal{L}_1$ and $\mathcal{L}_2$ change by $-2$ and $+2$, respectively.\\
\end{enumerate}

\begin{figure}[ht]
\scalebox{1.0}{\raisebox{-0pt}{$\vcenter{\hbox{\epsffile{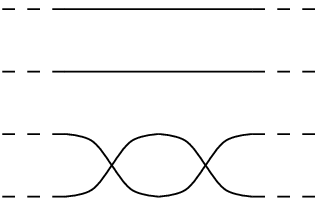}}}$}} \quad $\longrightarrow$ \quad
\scalebox{1.0}{\raisebox{-0pt}{$\vcenter{\hbox{\epsffile{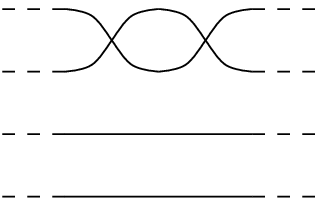}}}$}}
\caption{}
\end{figure}

These two operations allow us to remove at least half the number of self-crossings of the component $\mathcal{L}_2$. Hence we can realise all integer points up to the midpoint on the diagonal line segment between the points $(-1,\overline{tb}(L)-2lk(L)+1)$, $(\overline{tb}(L)-2lk(L)+1,-1)$ of the Thurston-Bennequin polytope of the link $L$. We conclude by remarking that every two-bridge link admits a symmetry that interchanges the two link components.  
\end{proof}

It is natural to ask whether the Thurston-Bennequin polytope of a two-component link can always be described by three linear inequalities, and whether all its integer points are realised by Legendrian framings. We do not know the answer to these questions. However, there do exist links whose Thurston-Bennequin polytope has an even simpler shape: positive braid links. According to Bennequin~\cite{Be} (see also~\cite{Ta}), every positive braid link has a natural Legendrian representative whose Thurston-Bennequin number is maximal and equals the number of crossings minus the number of strands of the braid. Moreover, the Thurston-Bennequin number of every single component of that representative is maximal. Thus the Thurston-Bennequin polytope of a two-component positive braid link is described by two inequalities, only:
$$x_1 \leq \overline{tb}(L_1), x_2 \leq \overline{tb}(L_2).$$
The same is evidently true for all two-component split links.

\section{Application related to quasipositivity}

Originally, a surface embedded in $\R^3$ was called quasipositive, if it was isotopic to a certain positively braided surface. Rudolph's characterisation results say that quasipositivity of surfaces is inherited by incompressible sub-surfaces and is preserved under the plumbing operation (\cite{Ru1}, \cite{Ru2}). Moreover, every quasipositive surface can be embedded incompressibly on the fibre surface of a torus link. The latter are plumbings of finitely many positive Hopf bands. The following definition is therefore justified.

\begin{definition} A surface embedded in $\R^3$ is called quasipositive, if it is obtained from a disc by plumbing finitely many positive Hopf bands and passing to an incompressible sub-surface.
\end{definition} 

An alternative definition for quasipositivity was recently formulated in the setting of contact geometry \cite{BI}: quasipositive surfaces are precisely the Legendrian ribbons in the standard contact $\R^3$. For our purposes, it suffices to know that every Legendrian link defines a Legendrian ribbon, by replacing its components by disjoint annuli with the obvious framings. 

\begin{proof}[Proof of Theorem~2]
We will construct an embedded surface, starting from the Legendrian ribbon associated to the front projection shown in Figure~6. This front realises the maximal Thurston-Bennequin number of the corresponding link (the negative torus link of type $(2,-4)$), by Ng \cite{Ng}. Therefore, if we remove the full twist encircled in Figure~6, we obtain a surface, say $F$, which is not a Legendrian ribbon, hence a non-quasipositive surface. Nevertheless, both incompressible sub-annuli of $F$ are obviously quasipositive.

\begin{figure}[ht]
\scalebox{0.8}{\raisebox{-0pt}{$\vcenter{\hbox{\epsffile{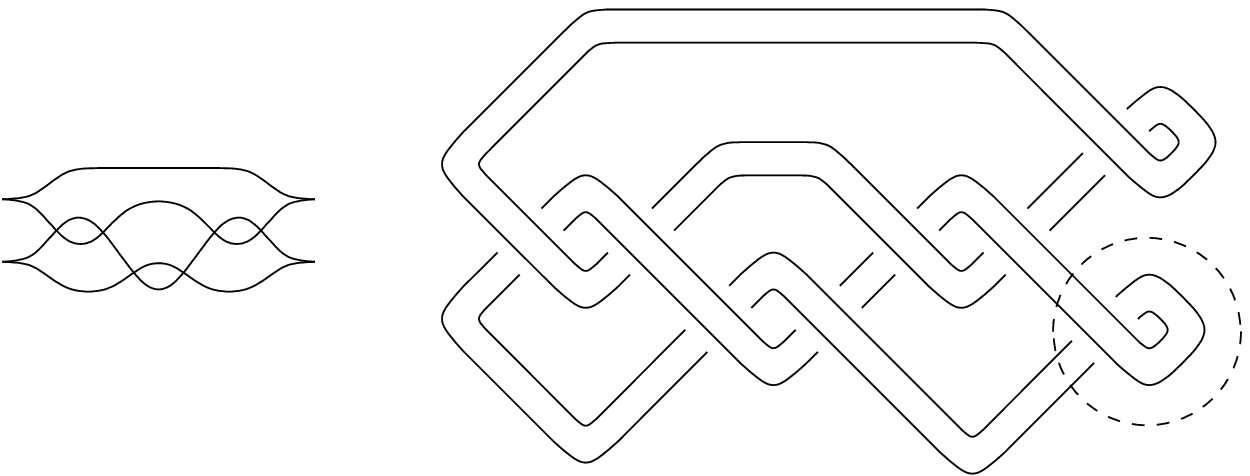}}}$}} 
\caption{}
\end{figure}

Connecting these two sub-annuli by a small band, as shown at the top of Figure~7, we obtain a connected surface $F'$, which is again non-quasipositive, since it contains $F$ as an incompressible sub-surface. We claim that all incompressible sub-annuli of $F'$ are quasipositive. 
Every incompressible sub-annulus of $F'$ is parallel to the boundary of $F'$. Indeed, the topological type of the surface $F'$ is a pair of pants, thus the core curve of an incompressible sub-annulus divides $F'$ into a pair of pants and an annulus. We already checked the quasipositivity for two sub-annuli of $F'$. The third sub-annulus is shown at the bottom of Figure~7. It is manifestly a Legendrian ribbon, hence quasipositive. 
\end{proof}

The above construction works with other two-bridge links, as well, for example with all negative torus links of type $L(2,-2n)$, $n \geq 2$.

\begin{figure}[ht]
\scalebox{0.8}{\raisebox{-0pt}{$\vcenter{\hbox{\epsffile{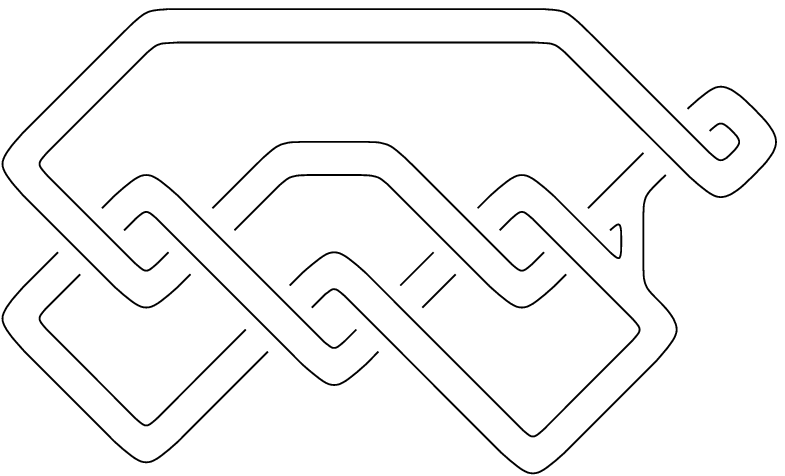}}}$}} 

\bigskip
\bigskip
\scalebox{0.8}{\raisebox{-0pt}{$\vcenter{\hbox{\epsffile{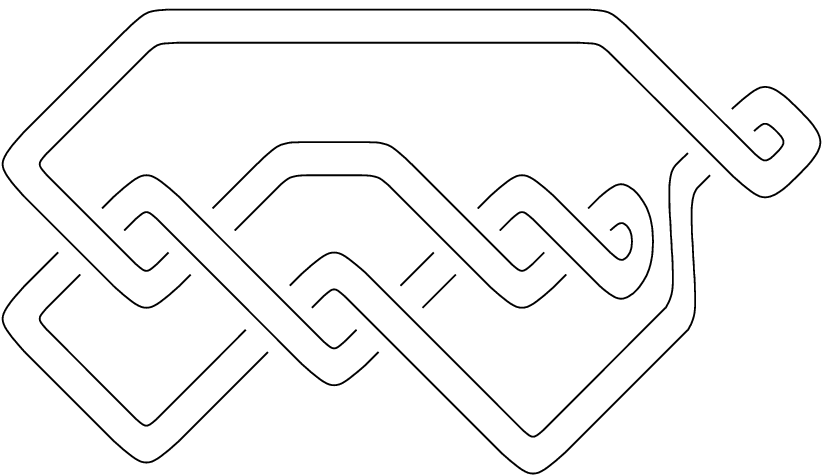}}}$}} 
\caption{}
\end{figure}

\bigskip
\noindent
Departement Mathematik, ETH Z\"urich, R\"amistrasse 101, CH-8092 Z\"urich, Switzerland.

\noindent
\texttt{sebastian.baader@math.ethz.ch}

\bigskip
\noindent
Mathematical Institute, Tohoku University, Sendai, 980-8578, Japan.

\noindent
\texttt{ishikawa@math.tohoku.ac.jp}

\end{document}